\documentclass[12 pt]{article}
\usepackage{amssymb,amsmath,amsfonts,amsthm}
\usepackage{xcolor}
\newcommand\BBP{{\mathbb {P}}}

\newcommand\BBR{{\mathbb {R}}}

\newcommand\BBE{{\mathbb {E}}}

\newcommand\E{{\mathbb {E}}}

\newtheorem{theorem}{Theorem}

\newtheorem{corollary}[theorem]{Corollary}

\newtheorem{lemma}[theorem]{Lemma}

\newtheorem{proposition}[theorem]{Proposition}
\newtheorem{remark}[theorem]{Remark}

\begin{document}


\title{On Berry Esseen type estimates for randomized Martingales in the non stationary setting}

\author{J. Dedecker\footnote{J\'er\^ome Dedecker, Universit\'e Paris Cit\'e, CNRS, MAP5, UMR 8145,
45 rue des  Saints-P\`eres,
F-75006 Paris, France.}, F. Merlev\`ede \footnote{Florence Merlev\`ede, LAMA,  Univ Gustave Eiffel, Univ Paris Est Cr\'eteil, UMR 8050 CNRS,  \  F-77454 Marne-La-Vall\'ee, France.}, 
M. Peligrad \footnote{Magda Peligrad, Department of Mathematical Sciences, University of Cincinnati, PO Box 210025, Cincinnati, Oh 45221-0025, USA.} and Vishakha  \footnote{Vishakha, Department of Mathematical Sciences, University of Cincinnati, PO Box 210025, Cincinnati, Oh 45221-0025, USA.}}

\maketitle

\begin{abstract} 
In this paper, we consider partial sums of triangular martingale differences weighted
by random variables drawn uniformly on the sphere, and globally independent of the
martingale differences. Starting from the so-called principle of conditioning and using some arguments developed by Klartag-Sodin and  Bobkov-Chistyakov-G\"otze, we give some upper bounds for the Kolmogorov distance between the distribution of these weighted sums and a Normal  distribution. Under some conditions on the conditional variances of the martingale differences, the obtained rates are always faster than those obtained  in case of usual partial
sums. 
\end{abstract}



\bigskip

\section{Introduction}
Let $(\Omega ,\mathcal{F},{\mathbb{P}})$ be a probability space and $(d_{j,n})_{1\leq j\leq n, n \geq 1}$ be a triangular array of martingale differences in $ {\mathbb L}^2(\Omega ,\mathcal{F},{\mathbb{P}})$, adapted to a triangular array of filtrations  $({\mathcal F}_{j,n})_{1\leq j\leq n}$ (so ${\mathcal F}_{j,n} \subset {\mathcal F}_{j+1,n}$ and 
$\E ( d_{j,n} | {\mathcal F}_{j-1,n}) =0$) and such that
\begin{equation}
\sum_{k=1}^{n}\BBE (d_{k,n}^{2})=n  \, . \label{Bobkov2}%
\end{equation}

In this paper, we are interested in the rates of convergence in the CLT for randomized weighted sums associated with  $(d_{j,n})_{1\leq j\leq n}$ in the following sense. Let 
$\theta =(\theta _{1},\dots ,\theta _{n})^t$ be a random vector defined on $(\Omega ,\mathcal{F},{\mathbb{P}})$, independent of $({\mathcal F}_{j,n})_{1\leq j\leq n, n \geq 1}$ and with  uniform distribution $\mu _{n-1}$ on
the unit sphere $S^{n-1}$ of ${\mathbb{R}}^{n}$ ($n\geq 2$).
Let $S_{n}(\theta)=\sum_{j=1}^{n}\theta_{j}d_{j,n}$. By independence,  taking the conditional expectation of $f(\theta, (d_{j,n})_{1\leq j\leq n})$ with respect to $\theta$ means keeping $\theta$ fixed and integrating with respect to the distribution of $(d_{j,n})_{1\leq j\leq n}$. We will designate by ${\mathbb P}_{|\theta}$, the conditional probability given $\theta$.  When we use conditional expectations, the equalities and inequalities are understood in the almost sure sense.

We shall give convergence rates in terms of the Kolmogorov distance, and more precisely for the quantity 
\[
 \BBE ( \kappa_{\theta}(P_{S_{n}(\theta)},P_{N_{\theta}}))  \quad \text{with} \quad \kappa_{\theta}(P_{S_{n}(\theta)},P_{N_{\theta}})=\sup_{t\in{\mathbb{R}}}%
|{\mathbb P}_{|\theta}(S_{n}(\theta)\leq t)-{\mathbb P}_{|\theta}(N_{\theta}\leq t)|\, ,
\]
where, conditionally to $\theta$,  the random variable $N_\theta$ is Gaussian with mean 0  and variance $\sum_{k=1}^{n} \theta_k^2\BBE (d_{k,n}^{2})$.  We will also denote by $N$ a centered 
Gaussian random variable with variance one. 

When  $(d_{k})_{k\in {\mathbb{Z}}}$ forms a sequence of independent centered
random variables in ${\mathbb{L}}^{4}$ with $\BBE (d_{k}^{2})=1$,  Corollary 3.2
in Klartag and Sodin \cite{KS} states that 
\begin{equation}
{\mathbb{E}}\Big (\kappa _{\theta }(P_{S_{n}(\theta )},P_{N})\Big )\leq 
\frac{c}{n^2}\sum_{k=1}^{n}{\mathbb{E}}(d_{k}^{4})\,.  \label{resultindependent}
\end{equation}
Therefore, in the iid setting, the upper bound in this Berry-Esseen type estimate is of order $n^{-1}$ whereas,  the classical Berry-Esseen theorem (meaning that we take 
$\theta_i = 1/ \sqrt{n}$ for any integer $i $ in $[1,n]$), gives a rate of order $n^{-1/2}$.  As quoted in Section 3 of \cite{DMP24} this type of randomized partial sums appears naturally  when we study the ordinary least square estimator of the slope in
a linear regression model with Gaussian design. As a consequence a super fast rate can be achieved 
for the Kolmogorov distance between the centered and normalized  least square estimator and its Gaussian limit, when the random design is Gaussian and the errors are iid and in ${\mathbb L}^4$.
In fact,  this also holds when the errors form a stationary sequence of  martingale differences satisfying some weak dependence conditions, see \cite[Corollary 2.3]{DMP24}.  

This type of phenomena concerning the Berry-Essen bound for randomized partial sums has  also been  studied in a series of papers by Bobkov, Chistyakov and G\"otze   \cite{BCG18,BCG20a,BCG20b} where the independence assumption has been relaxed. Let us describe the result given in Theorem 1.1 of \cite{BCG20b}. Let $(d_{k})_{k\in {\mathbb{Z}}%
}$ be a sequence of uncorrelated centered random variables with variance
one, satisfying the following second order correlation condition: there exists a constant $\Lambda $ such that, for any $n \geq 1$ and   any collection $%
a_{ij}\in {\mathbb{R}}$, 
\begin{equation}
\mathrm{Var}\Big (\sum_{i,j=1}^{n}a_{ij}d_{i}d_{j}\Big )\leq \Lambda
\sum_{i,j=1}^{n}a_{ij}^{2}\, . \label{condsecondorder}
\end{equation}%
If moreover the random vector  $(d_{1},\ldots ,d_{n})$
has a symmetric distribution then 
\begin{equation}
{\mathbb{E}}\Big (\kappa _{\theta }(P_{S_{n}(\theta )},P_{N})\Big )\leq 
\frac{c\log n}{n}\Lambda \,.  \label{resultsecondorder}
\end{equation}
The main restriction in this result is probably the fact that the
distribution of $(d_{1},\ldots ,d_{n})$ is assumed to be symmetric. In \cite[%
Chapter 17.4]{BCGbook} the authors consider the case of non-symmetric
distributions, but this leads to an additional term in the upper bound \eqref{resultsecondorder} (see \cite[Proposition 17.4.1]{BCGbook}). 
Let us comment now on the second order condition \eqref{condsecondorder}. As shown in \cite{BCG20b},  it can be verified
for random vectors satisfying a Poincar\'{e}-type inequality with
positive constant (see \cite[Proposition 3.4]{BCG20b}). On another hand, as quoted in \cite{DMP24}, if $(d_{k})_{k\in {\mathbb{Z}}}$ is a sequence of martingale
differences such that $\sup_{i\geq 1}{\mathbb{E}}(d_{i}^{4})<\infty $, then condition \eqref{condsecondorder} is satisfied under some natural dependence conditions.  

In the recent paper \cite{DMP24},  Dedecker, Merlev\`ede and Peligrad have provided a new method  allowing  to show that stationary sequences $%
(d_{k})_{k\in {\mathbb{Z}}}$ of martingale differences  satisfy an upper
bound of the type \eqref{resultindependent} (up to some logarithmic terms) without requiring that the law
of the vector $(d_{1},\ldots ,d_{n})$ is symmetric, provided that the sequence $(d_{k})_{k\in {\mathbb{Z}}}$ satisfies a weak dependence condition.  Their proof is based on Lindeberg's method combined with a variant of Berry-Esseen's smoothing inequality, which is adapted to the type of randomisation we use (see Lemma 5.2 in 
\cite{BCG18}). Note, however, that the arguments they developed for Lindeberg's  method to work require that the  sequence $(d_{k})_{k\in {\mathbb{Z}}}$ is   close to a stationary regime.

In this paper we are interested  in the extension of the result of Klartag-Sodin (up to some additional logarithmic terms) to triangular array of martingale differences $(d_{i,n})_{ 1\leq i \leq n}$ that are non necessarily stationary.  As a starting point, we shall use the  Principle of Conditioning saying heuristically that starting with a
limit theorem for independent random variables, if one replaces the
expectations by conditional expectations with respect to the past, summation
to constants by summation to stopping times, and the convergence of numbers by
convergence in probability, one obtains a correct statement in the dependent
setting. This heuristic was expressed as a rigorous theorem by Jakubowski
\cite{Ja86}, and appears for instance in the books by De la Pe\~{n}a and Gin\'{e} \cite[Chap. 7.1]{DG99} and Merlev\`{e}de et al. \cite[Section 2.3]{MPU19}. This principle was used to rigorously transport several limit theorems such as the central limit theorem or the moderate deviation principle, from the independent setting to the martingale one.

The paper is organized as follows. In Section \ref{SectionMR} we present our main 
result concerning the approximation of  ${ S}_n(\theta)$ by a Gaussian random variable in term of the Kolmogorov distance,  in the case where $(d_{i,n})_{1\leq i \leq n}$ is a triangular array of (non necessarily stationary)  martingale differences. In Subsection 2.1 we consider an example of a nonstationary sequence of martingale differences given in \cite{DMR22}, for which the rates in the usual Berry-Esseen theorem is at best of order $n^{-1/3}$, while the rate for randomized sums is at most of order $n^{-1}(\log n)^2$. In Section \ref{Secexample}, we give an example of a nonstationary ARCH model to which our results apply. The proofs are postponed to Section \ref{sectionproofs}. 

We shall use most of the times the notation $\ll$ instead of the Vinogradov
symbol $O$. Note that all our results hold if we replace the random vector $\theta=(\theta
_{1},\dots,\theta_{n})$ by ${\tilde \xi}= \Vert \xi \Vert_e^{-1} \xi$ where $\xi = (\xi_{1},\dots,\xi_{n})$ with $(\xi_i)_{1 \leq i \leq n}$  iid centered and standard Gaussian r.v.'s independent of $(d_{j,n})_{1\leq j\leq n}$ and $\Vert \xi \Vert_e^2 = \sum_{i=1}^n \xi_i^2$ stands for the Euclidean norm. Indeed, it is well-known that ${\tilde \xi}$  has uniform distribution on the unit sphere $S^{n-1}$ of
${\mathbb{R}}^{n}$.

\section{Main results} \label{SectionMR}



Recall that $S_{n}(\theta)=\sum_{j=1}^{n}\theta_{j}d_{j,n}$, where $\theta=(\theta
_{1},\dots,\theta_{n})$ is independent of $(d_{j,n})_{1\leq j\leq n}$
 and has uniform distribution on the unit sphere $S^{n-1}$ of
${\mathbb{R}}^{n}$ ($n\geq2$).  To soothe the notation we shall denote $d_{j}=d_{j,n}$  and ${\mathcal F}_{j}={\mathcal F}_{j,n}$. 
Let ${\mathbb E}_{j}$ be the conditional expectation with respect to ${\mathcal F}_j$. We introduce also the following notations: 
\begin{equation}
\sup_{1 \leq k \leq n }{\mathbb E}(d_{k}^{4}) = \alpha_n \ \text{ and } \  {\mathbb E} \left(  \frac{1}{n}%
\sum_{j=1}^{n}{\mathbb E}^2_{j-1} ( d_{j}^{3} ) \right)^{1/2} = \beta_n  \, , \label{Bobkov1}
\end{equation}
\begin{equation}
\sigma_{4}^{2} (n) =\frac{1}{n}\mathrm{var}\left (\sum_{k=1}^{n}d_{k}^{2}\right )  \, , 
\label{Bobkov3}%
\end{equation}
and 
\begin{equation} \label{condcondexp}
\sum_{ j = 1}^n \Vert {\mathbb E}_{j-1}(d_{j}^{2}) - {\mathbb E}(d_{j}^2)  \Vert_1 =  \gamma_n \, .
\end{equation}

Our main result is the following:

\begin{theorem}
\label{Th martingale}Assume that $(d_{i})_{1\leq i\leq n}$ is a vector of
martingale differences in ${\mathbb L}^4$, adapted to an array $(\mathcal{F}_{i})_{0\leq i\leq n}$
of increasing sigma fields and satisfying the condition \eqref{Bobkov2}. Then%
\[
 \BBE ( \kappa_{\theta}(P_{S_{n}(\theta)},P_{N_{\theta}})) 
\ll
\frac{  1+ v_n }{n},
\]
where the conditional distribution of  $N_\theta$ given $\theta$ is a normal distribution with mean 0 and variance $\sum_{i=1}^n \theta_i^2 {\mathbb E} (d_i^2)$ and 
\[
v_n =  \gamma_n (\log n )+  \beta_n (\log n)^{3/2}+ \alpha_n(\log n)^{2}  + \sigma_{4}^{2} (n) (\log n) \, .
\]
\end{theorem}

\begin{remark} \label{remarkonGaussian} Let $N$ be a standard centered Gaussian random variable.  If we are interested in giving an upper bound for the quantity 
\[
{\mathbb E}\sup_{t\in{\mathbb{R}}}|{\mathbb P}_{|\theta}(S_{n}(\theta)\leq t)-{\mathbb P}(N \leq t)| \, , 
\]
we can use Theorem \ref{Th martingale} together with the following upper bound whose proof is given in Subsection \ref{Sectionremark2}:  for any $n \geq 2$, 
\begin{equation} \label{boundGaussian}
{\mathbb E}\sup_{t\in{\mathbb{R}}}|{\mathbb P}_{|\theta}(N_\theta\leq t)-{\mathbb P}(N \leq t)|  \leq \frac{35 \sqrt{5}}{4}  \frac{    \Big ( \sum_{k=1}^n ( b_{k,n}^2 -1)^2 \Big )^{1/2}  }{n-1} \, , 
\end{equation}
where $b_{k,n}^2 = {\mathbb E} (d_{k,n}^2)$. Now, if we are interested in  
$
\sup_{t\in{\mathbb{R}}}|{\mathbb P}(S_{n}(\theta)\leq t)-{\mathbb P}(N \leq t)| 
$, note that
\begin{equation} \label{boundGaussian2}
\sup_{t\in{\mathbb{R}}}| {\mathbb E} ({\mathbb P}_{|\theta}(N_\theta\leq t))-{\mathbb P}(N \leq t)|  \ll \frac{  1+  \alpha_n ( \log n)^2  }{n} \, . 
\end{equation}
\end{remark}
In Theorem \ref{Th martingale} the upper bound makes use of the quantity $\gamma_n$ defined in \eqref{condcondexp},  and to get an interesting upper bound, we need to have ${\mathbb E}_{j-1}(d_{j}^{2})$ ``not too far'' in ${\mathbb L}^1$ from $  {\mathbb E}(d_{j}^2) $. Alternatively, we can  introduce the following additional notation: 
\begin{equation}
c^2_2(n ) =\frac{1}{n} \sum_{k=1}^{n} \sum_{\ell=1}^{n} | \mathrm{Cov}(d_{k}^{2} , d^2_{\ell})  |  \, ,
\label{Bobkov3alter}%
\end{equation}
and ask this quantity to be ``not too big".  Requiring  $c^2_2(n ) $ to be  ``not too big" is usually less stringent than imposing $\gamma_n$ to be small because 
$c^2_2(n )$ only requires fast enough decorrelation between $d_{k}^{2}$ and $  d^2_{\ell} $ when $|k - \ell |$ is large enough.  Note also that $\sigma_{4}^{2} (n)  \leq c^2_2(n) $.  
Using the quantity $c^2_2(n )$ rather than $\gamma_n$ gives the following Berry-Esseen bound: 
\begin{theorem}
\label{Th martingale2}Assume that $(d_{i})_{1\leq i\leq n}$ is a vector of
martingale differences in ${\mathbb L}^4$, adapted to an array $(\mathcal{F}_{i})_{0\leq i\leq n}$
of increasing sigma fields and satisfying the condition \eqref{Bobkov2}. Then%
\[
{\mathbb E}\sup_{t\in{\mathbb{R}}}|{\mathbb P}_{|\theta}(S_{n}(\theta)\leq t)-{\mathbb P}(N \leq t)|\ll
\frac{  1+ v'_n }{\sqrt{n}},
\]
where $N$ is a centered  standard normal variable and 
\[
v'_n = ( c_2(n ) \sigma_{4}^{2} (n) )\log n +  \beta_n (\log n)^{3/2}+ \alpha_n(\log n)^{2}  
\, .
\]
\end{theorem}
As we can see,   if the quantities $ c_2(n ) $, $\beta_n$ and $ \sigma_{4}^{2} (n)$ are uniformly bounded in $n$ by a finite constant, the bound in the Berry-Esseen theorem for this type of randomized martingales is of order $(\log n)^{2} / \sqrt{n}$. This rate is slower than the one we could obtain by applying  Theorem \ref{Th martingale} but as a counterpart asking the quantity $ c_2(n ) $ to be  uniformly bounded in $n$ is less stringent than making a similar assumption on $\gamma_n$. 

\subsection{Comparison with rates in the usual Berry-Essen bound for non stationary Martingales.} 

For a non stationary sequence of martingale differences, even in case of constant conditional variances and uniformly bounded third moments, the upper bound in
the classical Berry-Esseen inequality cannot be better than $n^{-1/4}$ (see \cite{Bo82}). In the same context, if there is a uniformly bounded moment of order $p>2$, Dedecker, Merlev\`ede and Rio 
\cite{DMR22} proved the following lower bound.  

\begin{proposition}[\cite{DMR22}] \label{CE} Let $p >2$ and $n \geq 20$. There exists $(d_1, \ldots, d_n)$ such that 
\begin{enumerate}
\item $\E(d_k | \sigma( d_1, \dots , d_{k-1} ) )=0 $ and $\E(d^2_k | \sigma( d_1, \dots , d_{k-1} )) =1  $ a.s.,
\item $ \sup_{1 \leq k \leq n} \E(|d_k|^p ) \leq \E ( |Y|^p) + 5^{p-2} $ where $Y \sim {\mathcal N} (0,1) $, 
\item $ \sup_{t\in \BBR}  \big| \BBP ( S_n \leq t \sqrt{n} ) - \Phi (t) \big| \geq 0.06 \,\, n^{-(p-2)/(2p-2)}$, where $S_n= \sum_{k=1}^n d_k$. 
\end{enumerate}
\end{proposition}

In particular, if we take $p=4$ in the above proposition, the lower bound in the classical Berry-Essen inequality is  of order $n^{-1/3}$. 
Now, for the vector  $(d_1, \ldots, d_n)$ constructed in the proof of the above proposition, we have that $\gamma_n =0$ and the other quantities ($\alpha_n$, $\beta_n$ and $\sigma_4^2 (n)$) are uniformly bounded. In addition, since $\sum_{i=1}^n \theta_i^2 \E (d_i^2) =1$, it follows that 
\[
\E\sup_{t\in{\mathbb{R}}}|{\mathbb P}_{|\theta}(S_{n}(\theta)\leq t)-{\mathbb P}(N \leq t)|\ll
\frac{  ( \log n)^2 }{n} \, .
\]   

\subsection{Example}  \label{Secexample} Let us consider $(\varepsilon_k)_{k \geq 1}$ a sequence of independent r.v.'s such that 
${\mathbb P} ( \varepsilon_k =0) = 1 - 2p_k$ and  ${\mathbb P} ( \varepsilon_k =a_k) = {\mathbb P} ( \varepsilon_k =-a_k) = p_k$ where 
$(a_k)_{k \geq 1}$ and $(p_k)_{k \geq 1}$ are sequences of positive reals such that $2a_k^2 p_k =1$, $2 p_k \leq 1$ and $\sum_{k \geq 1} p_k =1$. For any $i \geq 1$, we consider the ARCH model 
\[
d_i = \eta_{i-1} \varepsilon_i \, , 
\]
where $\eta_0=1$ and, for any $k \geq 1$,
\[
\eta_k^2 = \alpha_k^{-1} \sum_{\ell=1}^{k} u_\ell \varepsilon^2_{k +1 -\ell} \, , 
\]
with $(u_\ell)_{\ell \geq 1}$ a sequence of positive reals and $\alpha_k = \sum_{\ell=1}^{k} u_\ell$.  

We note that $(d_k)_{k \geq 1}$ is a sequence of martingale differences adapted to the filtration $({\mathcal F}_k )_{k \geq 1}$ where 
$\mathcal F_k  = \sigma ( \varepsilon_i, 1 \leq i \leq k)$.  Moreover,  since ${\mathbb E} (\varepsilon_i^2) =1$, we have ${\mathbb E} (d_i^2) =1$. In addition, for any $k \geq 1$ and any 
$\ell \geq k$,
\begin{multline*}
{\mathbb E}_k ( d_\ell^2) -1 = \alpha_{\ell-1}^{-1} \sum_{j=1}^{\ell-1} u_j {\mathbb E}_k ( \varepsilon^2_{\ell -j} \varepsilon^2_{\ell}  ) -1 \\
=   \alpha_{\ell-1}^{-1} \sum_{j=1}^{\ell-k-1} u_j {\mathbb E} ( \varepsilon^2_{\ell -j} \varepsilon^2_{\ell}  ) + \alpha_{\ell-1}^{-1} \sum_{j=\ell-k}^{\ell -1} u_j \varepsilon^2_{\ell -j} {\mathbb E}( \varepsilon^2_{\ell} )   -1  
=  \alpha_{\ell-1}^{-1} \sum_{j=\ell-k}^{\ell -1} u_j  \big ( \varepsilon^2_{\ell -j}   -1  \big ) \, .
\end{multline*}
Using the fact that $ a_k^2 = 1/(2 p_k ) \geq 1$ for any $k \geq 1$, we get 
\begin{multline} \label{petitcalcul}
\Vert {\mathbb E}_k ( d_\ell^2) -1 \Vert_1   \leq   \alpha_{\ell-1}^{-1} \sum_{j=\ell-k}^{\ell -1} u_j   \{ (1-2p_{\ell -j}) + 2 (a^2_{\ell -j} -1) p_{\ell -j} \}   
  =  2 \alpha_{\ell-1}^{-1}  \sum_{j=\ell-k}^{\ell -1} u_j  (1-2p_{\ell -j}) \, .
\end{multline}
We take now, for any $i \geq 1$,
\[
p_i = \frac{1}{2}  \Big ( 1 -  \frac{1}{i+1} \Big )  \mbox{ and } u_i = (i+1)^{-\gamma} \mbox{with }\gamma >2 \, .
\]
It follows from \eqref{petitcalcul} that, for any $\ell \geq 2$, 
\[
\Vert {\mathbb E}_{\ell-1} ( d_\ell^2) -1 \Vert_1  \ll {\ell}^{-1} \, , 
\]
proving that 
\[
\gamma_n = n^{-1}\sum_{\ell=1}^n \Vert {\mathbb E}_{\ell-1} ( d_\ell^2) -1 \Vert_1 \leq (\log n)/n \, .
\]
Note also that since, for any $k \geq 1$, $2a_k^2 p_k =1$, by the selection of $p_k$, we have that $\sup_{k \geq 1} a_k^2 \leq 2$ implying that 
$\sup_{k \geq 1}  \Vert d_k \Vert_{\infty} \leq 2$. It follows that the sequences $(\alpha_n)_{n \geq 1}$ and $(\beta_n)_{n \geq 1}$ defined in \eqref{Bobkov1}  satisfy $\alpha_n = 2^4$ and $\beta_n= 2^3$. It remains to give an upper bound for $\sigma_4^2(n)$.  We have, by \eqref{petitcalcul}, 
\begin{multline*}
\sigma_4^2(n) \leq \frac{2}{n} \sum_{k=1}^n   \sum_{\ell=k}^n \Vert d^2_k {\mathbb E}_k ( d_\ell^2 -1)   \Vert_1 
\leq \frac{4}{n} \sum_{k=1}^n   \sum_{\ell=k}^n \Vert  {\mathbb E}_k ( d_\ell^2 -1)   \Vert_1 \\
\leq  \frac{8}{n} \sum_{k=1}^n   \sum_{\ell=k}^n  \alpha_{\ell-1}^{-1}  \sum_{j=\ell-k}^{\ell -1} \frac{1}{(j+1)^{\gamma}}  \frac{1}{ \ell-j+1} \leq C \, , 
\end{multline*}
where $C$ does not depend on $n$.  Hence, applying Theorem \ref{Th martingale}, we derive that 

\[
{\mathbb E}\sup_{t\in{\mathbb{R}}}|{\mathbb P}_{|\theta}(S_{n}(\theta)\leq t)-{\mathbb P}(N \leq t)|\ll
\frac{  (\log n)^2 }{n} \, .
\]

\section{Proofs} \label{sectionproofs}

\subsection{A preliminary result}

Since our method of proof is based on the Principle of Conditioning, we start by recalling some useful results in the independent setting. In what follows, $(Y_{j,n})_{1 \leq j \leq n, n \geq 1}$ is a triangular array of independent random variables that are centered and in ${\mathbb L}^4$. We shall use the following notations:%
\[
R_{2}^{2}=R_{2}^{2}(\theta,n)=\sum_{j=1}^{n}\theta_{j}^{2}{\mathbb E}(Y_{j,n}^{2}),
\]%
\[
R_{3}^{3}=R_{3}^{3}(\theta,n)=\sum_{j=1}^{n}\theta_{j}^{3}{\mathbb E}(Y_{j,n}^{3}),
\]%
\[
R_{4}^{4}=R_{4}^{4}(\theta,n)=\sum_{j=1}^{n}\theta_{j}^{4}{\mathbb E}(Y_{j,n}^{4}),
\]
and%
\[
f_{j}(s)={\mathbb E}(\mathrm{e}^{\mathrm{i}sY_{j,n}}).
\]

\begin{lemma}[Klartag-Sodin \cite{KS}]
\label{Taylor} For all $\theta$ and any  
$t\leq\min(R_{4}^{-1}( \theta, n),|R_{3} ( \theta , n) |^{-1})$ we have
\[
\Big \vert \prod_{j=1}^{n}f_{j}(\theta_{j}t)-\mathrm{e}^{-R_{2}%
^{2} ( \theta , n ) t^{2}/2}\Big \vert \leq  c \mathrm{e}^{-R_{2} ^{2} ( \theta ,n ) t^{2}/2}\left(
|R_{3}^{3}(\theta,n)|\cdot|t|^{3}+R_{4}^{4}(\theta, n )t^{4}\right)  \, ,
\]
where $c$ is a universal constant.
\end{lemma}
The proof follows, with small changes, the initial steps of the proof of Lemma 2.1 in \cite{KS} before integrating with $t$, and is therefore left to the reader. 

\subsection{The Principle of Conditioning for randomized partial sums}

The proof of Theorems \ref{Th martingale} and \ref{Th martingale2} is based on the principle of
conditioning going back to Jakubowski  \cite{Ja86} (see also \cite[Chap. 7.1]{DG99} and \cite[Section 2.3]{MPU19}). The
idea of the proof is to reduce the proof for the martingale case to the independent case.

Let $\mathcal{F}_{j, \theta}=\mathcal{F}_{j}\vee \sigma(\theta)$. Assuming that $\Omega$ is rich enough, there are a sigma algebra ${\mathcal G}$ and random variables
$(U_{j})_{1\leq j\leq n}$  such that: the conditional distribution of  $U_{j}$ given ${\mathcal G}$ is the conditional distribution of  $\theta_jd_{j}$
given $\mathcal{F}_{j-1, \theta}$, and $U_1, \ldots, U_n$ are conditionally independent given ${\mathcal G}$ (see for instance Proposition 6.1.5 in \cite{DG99} for the existence of such a sequence).
Denote by  $E_{{\mathcal G}}$ the
conditional expectation with respect to ${\mathcal G}$. Note that
\begin{equation}
f_{n,\theta}(t):=E_{{\mathcal G}}\left(  \mathrm{e}^{it\sum_{j=1}^{n}U_{j}%
}\right)  =\prod\nolimits_{j=1}^{n}E_{{\mathcal G}}\left(  \mathrm{e}^{itU_{j}%
}\right)  =\prod\nolimits_{j=1}^{n}{\mathbb E}(\mathrm{e}^{\mathrm{i}t\theta_{j}d_{j}%
}|\mathcal{F}_{j-1, \theta}) \label{defproduct}%
\end{equation}
(see for instance the last line of page 328 in \cite{DG99}.) 
For any martingale differences $\left(  d_{j}\right)  $ adapted to
$(\mathcal{F}_{j})$ \ denote
\begin{equation}
T_{2}^{2} = T_{2}^{2}(\theta ,n )=\sum_{j=1}^{n}\theta_{j}^{2}{\mathbb E}_{j-1}\left(  d_{j}^{2}\right)  ,\text{
}T_{3}^{3} = T_{3}^{3}(\theta ,n )=\sum_{j=1}^{n}\theta_{j}^{3}{\mathbb E}_{j-1}\left(  d_{j}^{3}\right)
 \label{def T}%
\end{equation}
and
\begin{equation}
T_{4}^{4} = T_{4}^{4}(\theta ,n )=\sum_{j=1}^{n}\theta_{j}^{4}{\mathbb E}_{j-1}\left(  d_{j}%
^{4}\right)  . \label{def Tbis}%
\end{equation}
Lemma \ref{Taylor} applied to the conditionally independent random variables
$(U_{j})_{1\leq j\leq n}$, leads to the following corollary.

\begin{corollary}
\label{Cor product}Assume that $(d_{i})_{1\leq i\leq n}$ is a vector of
martingale differences adapted to an array $(\mathcal{F}_{i})_{0\leq i\leq n}$
of increasing sigma fields. Then, for any   $t\leq\min(T_{4}^{-1},|T_{3}|^{-1})$, 
\[
\left\vert \prod_{j=1}^{n}{\mathbb E}(\mathrm{e}^{\mathrm{i}t\theta_j d_{j}%
}|\mathcal{F}_{j-1, \theta})-\mathrm{e}^{-T_{2}^{2}t^{2}/2}\right\vert \\
  \leq c\mathrm{e}^{-T_{2}^{2}t^{2}/2}\left(  |t|^{3}|T_{3}^{3}|+t^{4}%
T_{4}^{4}\right)  ,
\]
where $c>1$ is a universal constant.
\end{corollary}

\begin{remark}
\label{rem cor}As a consequence of Corollary \ref{Cor product} we can find
a universal constant $K\leq 1$ such that for $t\leq K\min(T_{4}^{-1},|T_{3}%
|^{-1}),$we have%
\begin{equation}
\prod_{j=1}^{n}|{\mathbb E}(\mathrm{e}^{\mathrm{i}t\theta_j d_{j}}|\mathcal{F}_{j-1, \theta}%
)|\geq\frac{1}{2}\mathrm{e}^{-T_{2}^{2}t^{2}/2}. \label{product larger}%
\end{equation}

\end{remark}

Now, if $\prod_{j=1}^{n}|{\mathbb E}(\mathrm{e}^{\mathrm{i}t\theta_j d_{j}}|\mathcal{F}%
_{j-1})|\neq0,$ we could use the well-known remarkable martingale property (it
is relation (6.24) in \cite{DG99}): 
\[
{\mathbb E}\left(  \frac{\mathrm{e}^{\mathrm{i}tS_{n}(\theta)}}{\prod_{j=1}%
^{n}{\mathbb E}(\mathrm{e}^{\mathrm{i}t\theta_j d_{j}}|\mathcal{F}_{j-1, \theta})}\Big | \theta \right)  =1,
\]
and use Corollary \ref{Cor product} to find Berry Esseen bounds for
${\mathbb E}(\mathrm{e}^{\mathrm{i}tS_{n}(\theta)}).$

To be able to rigorously use this approach, we shall construct an additional
vector of martingale $\left(  d_{j}^{\prime}\right)  $ differences adapted to
$(\mathcal{F}_{j})_{0\leq j\leq n}$ satisfying 

\textbf{ }
\begin{equation}
|f_{n,\theta}^{\prime}(t)|>\frac{1}{2}\mathrm{e}^{- t^{2} ( V_n^2 ( \theta) + t^{-2} )  /2
} = \frac{1}{2 \sqrt{\mathrm{e}}}\mathrm{e}^{- t^{2}  V_n^2 ( \theta)  /2
} \, , \label{removable}%
\end{equation}
where
\begin{equation}
f_{n,\theta}^{\prime}(t)= \prod_{j=1}^{n}{\mathbb E}(\mathrm{e}%
^{\mathrm{i}t\theta_{j}d_{j}^{\prime}}|\mathcal{F}_{j-1, \theta})  \mbox{ and } V_n^2 ( \theta) := \sum_{k=1}^n \theta_k^2 {\mathbb E} (d_k^2) .
\label{def product prime}%
\end{equation}
For this new sequence of martingale differences we find Berry-Esseen
bounds for $S_{n}^{\prime}(\theta)=\sum_{j=1}^{n}\theta_{j}d_{j}^{\prime}$ and
then compare $S_{n}^{\prime}(\theta)$ with $S_{n}(\theta)$.  Here are  the details.

\medskip

\noindent \textit{Construction of the auxiliary sequence $(d_{j}^{\prime})_{1\leq j \leq n}$ of martingale differences.} Define
\begin{equation}
A_{j}=A_{j}(\theta)= \Big \{ |f_{j,\theta}(t)|>\frac{1}{2 \sqrt{e}}\mathrm{e}^{- t^{2} V_n^2 ( \theta)   /2
} \Big \} \mbox{ and }d_{j}^{\prime}=d_{j}I_{A_{j}} . \label{def d prime}%
\end{equation}
Note that the $A_j$'s are ${\cal F}_{j-1, \theta}$ measurable and non-increasing in $j$ (for the inclusion). Also define for $(d^{\prime}_j)_{1\leq j \leq n}$ the quantities $T_{k}^{\prime}  = T_{k}^{\prime} ( \theta , n )$, $k=2,3,4$ as in
(\ref{def T}) and (\ref{def Tbis}).

This construction was used in \cite{Ja86,DG99,MPU19} where condition (\ref{removable}) was verified. For an
easy reference see relation (2.35) in \cite{MPU19}. Note that
(\ref{removable}) implies that%
\begin{equation} \label{removableprime}
|f_{k,\theta}^{\prime}(t)|>\frac{1}{2}\mathrm{e}^{- t^{2} ( V_n^2 ( \theta) + t^{-2} )  /2
} \text{ for all
}k\leq n.
\end{equation}
In addition, the following estimate will be needed in what follows:
\begin{equation} 
{\mathbb P}_{|\theta}(A_{n}^{c})\ll|t|^{3}{\mathbb E}_{|\theta}|T_{3}^{3}|+t^{4}{\mathbb E}_{|\theta}(T_{4}^{4}) + t^2  {\mathbb E}_{|\theta}\Big  | \sum_{k=1}^n   \theta_k^2 \big ( {\mathbb E}_{k-1} (d_k^2) - {\mathbb E} (d_k^2)   \big ) \Big  |  \, . \label{estimate Anc}%
\end{equation}
In case $t=0$, it is clear that ${\mathbb P}_{|\theta}(A_{n}^{c}) =0$, so that \eqref{estimate Anc} is true. Now 
to prove \eqref{estimate Anc} in case $t \neq 0$, we first write 
\[
{\mathbb P}_{|\theta}(A_{n}^{c}) \leq {\mathbb P}_{|\theta}(H_{n,K}^{c}) +  {\mathbb P}_{|\theta}\left(  |f_{n,\theta}(t)|<\frac{1}{2}\mathrm{e}^{-t^{2} (V_n^2 (\theta) +t^{-2})
/2},H_{n,K}(\theta)\right)  \, ,
\]
where 
\begin{equation}
H_{n,K}(\theta)= \big \{t\leq K \min(|T_{3}^{-1}|,T_{4}^{-1}) \big \} \, , \label{def H}%
\end{equation}
with $K$ defined in Remark \ref{rem cor}. By the Markov inequality,%
\begin{multline}
{\mathbb P}_{|\theta}(H_{n,K}^{c}(\theta))    \leq {\mathbb P}_{|\theta}(|t|^{3}>K^{3} |T_{3}|^{-3})+{\mathbb P}_{|\theta}(|t|^{4}>K^{4}T_{4}^{-4})  \\
  ={\mathbb P}_{|\theta}(T_{3}^{3}>K^{3}t^{-3})+{\mathbb P}_{|\theta}(T_{4}^{4}\geq K^{4}
t^{-4})\\
 \leq |t|^{3} K^{-3}{\mathbb E}_{|\theta}|T_{3}^{3}|+|t|^{4} K^{-4} {\mathbb E}_{|\theta}|T_{4}^{4}| \, . \label{evaluation H}
\end{multline}
Next noting that, on $H_{n,K}(\theta)$,   (\ref{product larger}) holds, we have, for any $t \neq 0$, 
\[
{\mathbb P}_{|\theta}\left(  |f_{n,\theta}(t)|<\frac{1}{2}\mathrm{e}^{-t^{2} (V_n^2 (\theta) +t^{-2})
/2},H_{n,K}(\theta)\right) \leq {\mathbb P}_{|\theta} (T_2^2 > V_n^2 (\theta) + t^{-2} ) \, .
\]
The estimate \eqref{estimate Anc} follows from the previous considerations and Markov's inequality. 

\medskip

In the next lemma, we analyse the characteristic function of $S_{n}(\theta)$. 

\begin{lemma}
\label{Lma1auxiliarymds} For any $t \in {\mathbb R}$, 
\begin{equation*}
\left\vert {\mathbb E}_{|\theta} (\mathrm{e}^{itS'_{n}(\theta)})-\mathrm{e}^{-t^{2} V_n^2 (\theta) /2}\right\vert
\ll B_n(\theta, t ) \text{ and }  \left\vert  {\mathbb E}_{|\theta} (\mathrm{e}^{itS_{n}(\theta)})-\mathrm{e}^{-t^{2} V_n^2 (\theta) /2}\right\vert
\ll B_n(\theta, t )\, ,
\end{equation*}
where 
\[
B_n(\theta, t ) =| t|^{3} {\mathbb E}_{|\theta} (|T_{3}^{3}|)+t^{4} {\mathbb E}_{|\theta} (|T_{4}^{4}|)   +  t^2 {\mathbb E}_{|\theta}\Big |\sum_{k=1}^n   \theta_k^2 \big (  {\mathbb E} _{k-1} (d_k^2) -  {\mathbb E} ( d_k^2 )\big ) \Big  | \, .
\]
\end{lemma}
\noindent {\bf Proof of Lemma \ref{Lma1auxiliarymds}.} Clearly we can assume that $t \neq 0$. First we  analyze ${\mathbb E}_{|\theta} \big ( \mathrm{e}^{itS_{n}%
^{\prime}(\theta)} \big )$. By (\ref{removable}), $f_{n,\theta}^{\prime
}(t)\neq0$. Hence we can start from a well-known identity, obtained by
recurrence, where we use the random variables $(d_{k}^{\prime}),$ namely:
\begin{equation}
{\mathbb E}_{|\theta}\left(  \frac{\mathrm{e}^{itS_{n}^{\prime}(\theta)}}{f_{n,\theta}^{\prime}(t)}\right)  =1.\label{identity}%
\end{equation}
Denote%
\[
\Delta_{n}(\theta)=\mathrm{e}^{itS_{n}^{\prime}(\theta)}-\mathrm{e}^{-t^{2} V_n^2 (\theta) /2}
\left(  \frac{\mathrm{e}^{itS_{n}^{\prime}(\theta)}}{f_{n,\theta}^{\prime}(t)}\right)   \, .
\]
Using (\ref{identity}), we evaluate%
\[
{\mathbb E}_{|\theta} \big ( \mathrm{e}^{itS_{n}^{\prime}(\theta)} \big )-\mathrm{e}^{-t^{2} V_n^2 (\theta) /2}={\mathbb E}_{|\theta} ( \Delta_{n}(\theta) ) \,  .
\]
Let
\[
H_{n,K}^{\prime}(\theta)=\{t\leq K\min(|T_{3}^{\prime-1}|,T_{4}^{\prime-1})\} \, ,
\]
with $K$ defined in Remark \ref{rem cor}.   So, by (\ref{product larger}),  on
$H_{n,K}^{\prime}(\theta)$, we have $|f_{n,\theta}^{\prime}(t)|\geq\frac{1}{2}\mathrm{e}^{- t^2T_{2}^{\prime2}%
(\theta,n)/2}$, which combined with \eqref{removableprime} implies that, on $H_{n,K}^{\prime}(\theta)$, 
\begin{equation}
|f_{n,\theta}^{\prime}(t)|\geq\frac{1}{2} \max \Big ( \mathrm{e}^{- t^2 T_{2}^{\prime2}%
(\theta,n)/2},  \mathrm{e}^{- t^2 ( V_n^2 (\theta) + t^{-2} ) /2 }
  \Big )  \, . \label{product larger prime}%
\end{equation}
With the above notations, we use the following decomposition:
\begin{equation}
{\mathbb E}_{|\theta} \big ( \mathrm{e}^{itS_{n}^{\prime}(\theta)}-\mathrm{e}^{-t^{2}V_n^2 ( \theta) /2} \big )={\mathbb E}_{|\theta} \big ( \Delta_{n}%
(\theta)I_{H_{n,K}^{\prime}(\theta)} \big ) +{\mathbb E}_{|\theta}\big ( \Delta_{n}(\theta)I_{H_{n,K}^{\prime
c}(\theta)} \big ) \, .\label{decopose}%
\end{equation}
By (\ref{removable}) we first notice that  
\begin{align*}
\big |{\mathbb E}_{|\theta} \big (  \Delta_{n}(\theta)I_{H_{n,K}^{\prime c}(\theta)} \big ) \big | & \leq {\mathbb E}_{|\theta}\left(  \left\{ |\mathrm{e}%
^{itS_{n}^{\prime}(\theta)}|+\mathrm{e}^{-t^{2}V_n^2 ( \theta) /2} \left\vert \frac
{\mathrm{e}^{itS_{n}^{\prime}(\theta)}}{f_{n,\theta}^{\prime}(t)}\right\vert
\right \}  I_{H_{n,K}^{\prime c}(\theta)} \right ) \\
& \leq (1 + 2 \sqrt{\mathrm{e}} ){\mathbb P}_{|\theta}(H_{n,K}^{\prime c}(\theta)) \, .
\end{align*}
Now
\[
{\mathbb P}_{|\theta}(H_{n,K}^{\prime c}(\theta))  \leq  {\mathbb P}_{|\theta}( H_{n,K}^{\prime c}(\theta) \cap A_n ) + {\mathbb P}_{|\theta} (A^c_n )\, .
\]
Note that on $A_n$, $H_{n,K}^{\prime c}(\theta) = H_{n,K}^{ c}(\theta) $ (since the $A_j$'s are non increasing in $j$). Therefore, taking into account  (\ref{evaluation H}) 
and  \eqref{estimate Anc} we get 
\begin{equation}
\big |{\mathbb E}_{|\theta} \big (  \Delta_{n}(\theta)I_{H_{n,K}^{\prime c}(\theta)} \big ) \big |\ll|t|^{3}{\mathbb E}_{|\theta}(|T_{3}^{3}|)+t^{4}{\mathbb E}_{|\theta}(T_{4}^{4}) + t^2  {\mathbb E}_{|\theta}\Big | \sum_{k=1}^n   \theta_k^2 \big ( {\mathbb E}_{k-1} (d_k^2) - {\mathbb E} (d_k^2)   \big ) \Big  | \, .\label{first part}%
\end{equation}
On another hand, note that
\begin{multline} \label{first part 2}
\big |{\mathbb E}_{|\theta} \big (  \Delta_{n}(\theta)I_{H_{n,K}^{\prime}(\theta)} \big ) \big |   \leq {\mathbb E}_{|\theta} \left ( \frac{\left\vert
f_{n,\theta}^{\prime}(t)-\mathrm{e}^{-t^{2}V_n^2 ( \theta )/2}\right\vert }{ \vert f_{n,\theta
}^{\prime}(t) \vert} I_{H_{n,K}^{\prime}(\theta) \cap A_n} \right ) \\ +  {\mathbb E}_{|\theta} \left ( \frac{\left\vert
f_{n,\theta}^{\prime}(t)-\mathrm{e}^{-t^{2}V_n^2 ( \theta )/2}\right\vert }{ \vert f_{n,\theta
}^{\prime}(t) \vert }I_{A^c_n} \right ) \, .
\end{multline}
By (\ref{removable}), 
\begin{equation}\label{truc1}
 {\mathbb E}_{|\theta} \left ( \frac{\left\vert
f_{n,\theta}^{\prime}(t)-\mathrm{e}^{-t^{2}V_n^2 ( \theta )/2}\right\vert }{ \vert f_{n,\theta
}^{\prime}(t) \vert }I_{A^c_n} \right )  \leq (1 + 2 \sqrt{\mathrm{e}} ){\mathbb P}_{|\theta}(A_n^{c}) \, .
\end{equation}
Now, on $H_{n,K}^{\prime}(\theta) $ we can use (\ref{product larger prime}) and
therefore,%
\begin{align*}
\big |{\mathbb E}_{|\theta} \big (  \Delta_{n}(\theta)I_{H_{n,K}^{\prime}(\theta)\cap A_n} \big ) \big | &  \leq 
{\mathbb E}_{|\theta} \left ( \frac{\left\vert
f_{n,\theta}^{\prime}(t)-\mathrm{e}^{-t^{2}V_n^2 ( \theta )/2}\right\vert }{ \vert f_{n,\theta
}^{\prime}(t) \vert }I_{H_{n,K}^{\prime}(\theta) \cap A_n} \right ) \\
&  \leq  2 {\mathbb E}_{|\theta} \left ( \mathrm{e}^{t^{2} \min ( T_{2}^{\prime2} (\theta , n) , V_n^2 ( \theta )  +t^{-2}) /2}\left \vert f_{n,\theta}^{\prime}%
(t)-\mathrm{e}^{-t^{2}V_n^2 ( \theta )/2}\right\vert I_{H_{n,K}^{\prime}(\theta)\cap A_n} \right ) \, .
\end{align*}
We continue the computations by introducing an intermediate term and use the
triangle inequality to obtain:%
\[
\big |{\mathbb E}_{|\theta}\big (  \Delta_{n}(\theta)I_{H_{n,K}^{\prime}(\theta)\cap A_n} \big ) \big | \leq 2 I_n(\theta) + 2 I\!\! I_n(\theta) \, , 
\]
where 
\[
 I_n(\theta)  = {\mathbb E} \left ( \mathrm{e}^{t^{2} T_{2}^{\prime2} (\theta , n) /2}\left \vert f_{n,\theta }^{\prime}%
(t)- \mathrm{e}^{-t^{2}T_{2}^{\prime2} (\theta , n)/2 }\right\vert I_{H_{n,K}^{\prime}(\theta)\cap A_n} \right )
\]
and
\[
 I\!\! I_n(\theta)  = {\mathbb E}_{|\theta} \left ( \mathrm{e}^{t^{2} \min ( T_{2}^{\prime2} (\theta , n) , V_n^2 ( \theta )  +t^{-2}) /2}\left | \mathrm{e}^{-t^{2}T_{2}^{\prime2} (\theta , n)/2 }-\mathrm{e}^{-t^{2}V_n^2 ( \theta )/2} \right \vert  I_{A_n} \right ) \, .
\]
By Corollary \ref{Cor product}  applied to $(d_{j}^{\prime})$, we get 
\begin{equation} \label{BoundIn} 
I_{n}(\theta)\ll {\mathbb E}_{|\theta} \big [(|t|^{3}|T_{3}^{{\prime}3}|+t^{4}T_{4}^{\prime4})  I_{A_n} \big ]   \ll  {\mathbb E}_{|\theta} \big [(|t|^{3}|T_{3}^{3}|+t^{4}T_{4}^{4}) \big ] \, ,
\end{equation}
since on $A_n$, $T_{3}^{{\prime}}= T_3$ and $T_{4}^{\prime}=T_4$. 

We give now an estimate for the quantity $ I\!\! I_n(\theta) $. With this aim, we introduce the following sets:
\[
B_1 (\theta) = \{   T_{2}^{\prime2} (\theta , n)  \leq  V_n^2 ( \theta )  \}  \, , \,  B_2 (\theta) = \{  V_n^2 ( \theta ) <  T_{2}^{\prime2} (\theta , n)  \leq   V_n^2 ( \theta ) + t^{-2}  \}
\]
and
\[
B_3 (\theta) = \{     T_{2}^{\prime2} (\theta , n) >  V_n^2 ( \theta ) + t^{-2}  \} \, .
\]
With these notations and setting 
\[
Q^{\prime}_n(\theta) =  \mathrm{e}^{t^{2} \min ( T_{2}^{\prime2} (\theta , n) , V_n^2 ( \theta )  +t^{-2}) /2}\left | \mathrm{e}^{-t^{2}T_{2}^{\prime2} (\theta , n)/2 }-\mathrm{e}^{-t^{2}V_n^2 ( \theta )/2}\right\vert  \, ,
\]
we have
\[
I\!\! I_n(\theta)  =\sum_{i=1}^3 {\mathbb E}_{|\theta} ( Q^{\prime}_n(\theta)  I_{B_i(\theta) \cap A_n} )  \, .
\]
Note first that 
\begin{multline*}
 {\mathbb E}_{|\theta} ( Q^{\prime}_n(\theta)  I_{B_1(\theta)\cap A_n} ) = {\mathbb E}_{|\theta} \left \{  \mathrm{e}^{t^{2} T_{2}^{\prime2} (\theta , n)  /2}\left ( \mathrm{e}^{-t^{2}T_{2}^{\prime2} (\theta , n)/2 }-\mathrm{e}^{-t^{2}V_n^2 ( \theta )/2}\right )   I_{B_1(\theta) \cap A_n} \right \}  \\
= {\mathbb E}_{|\theta}\left \{ \left ( 1-\mathrm{e}^{-t^{2} ( V_n^2 ( \theta ) - T_{2}^{\prime2} (\theta , n) /2}\right )   I_{A_n \cap B_1(\theta)} \right \} \, .
\end{multline*}
Now,  on $ B_1(\theta)$, $V_n^2 ( \theta ) - T_{2}^{\prime2} (\theta , n) \geq 0$ and on $A_n$, $ T_{2}^{\prime2} (\theta , n) = T_{2}^{2} (\theta) $. Hence, using the fact that, for $x \geq 0$, $1 - \mathrm{e}^{-x} \leq x$, it follows that 
\begin{equation}
{\mathbb E}_{|\theta} ( Q^{\prime}_n(\theta)  I_{B_1(\theta)\cap A_n} ) \ll  t^2 {\mathbb E}_{|\theta} \Big | \sum_{k=1}^n   \theta_k^2 \big ( {\mathbb E}_{k-1} (d_k^2) - {\mathbb E} (d_k^2)   \big ) \Big  |  \, . \label{estimate B1}%
\end{equation}
We give now an estimate of  $  {\mathbb E}_{|\theta} ( Q^{\prime}_n(\theta)  I_{B_2(\theta) \cap A_n} ) $.  We have 
\begin{multline*}
 {\mathbb E}_{|\theta} ( Q^{\prime}_n(\theta)  I_{B_2(\theta) \cap A_n } ) = {\mathbb E}_{|\theta} \left \{  \mathrm{e}^{t^{2} T_{2}^{\prime2} (\theta , n)  /2}\left ( \mathrm{e}^{-t^{2}V_n^2 ( \theta )/2} -  \mathrm{e}^{-t^{2}T_{2}^{\prime2} (\theta , n)/2 }\right )   I_{B_2(\theta) \cap A_n} \right \}  \\
 \leq  \sqrt{ \mathrm{e} }  {\mathbb E}_{|\theta} \left \{  \mathrm{e}^{t^{2} V_n^2 ( \theta )/2} \left ( \mathrm{e}^{-t^{2}V_n^2 ( \theta )/2} -  \mathrm{e}^{-t^{2}T_{2}^{\prime2} (\theta , n)/2 }\right )   I_{A_n \cap B_2(\theta)} \right \}  \, . \end{multline*}
Proceeding as above, we get that 
 \begin{equation}
 {\mathbb E}_{|\theta} ( Q^{\prime}_n(\theta)  I_{B_2(\theta) \cap A_n} )  \ll  t^2  {\mathbb E}_{|\theta} \Big | \sum_{k=1}^n   \theta_k^2 \big ( {\mathbb E}_{k-1} (d_k^2) - {\mathbb E} (d_k^2)   \big ) \Big  |  \, . \label{estimate B2}%
\end{equation}
It remains to deal with $  {\mathbb E}_{|\theta} ( Q^{\prime}_n(\theta)  I_{B_3(\theta) \cap A_n} ) $. We first write
\begin{multline*}
 {\mathbb E}_{|\theta} ( Q^{\prime}_n(\theta)  I_{B_3(\theta)\cap A_n} ) = {\mathbb E}_{|\theta} \left \{  \mathrm{e}^{t^{2}  ( V_n^2 ( \theta ) + t^{-2})/2  }\left ( \mathrm{e}^{-t^{2}V_n^2 ( \theta )/2} -  \mathrm{e}^{-t^{2}T_{2}^{\prime2} (\theta , n)/2 }\right )   I_{B_3(\theta)\cap A_n} \right \}  \\
 \leq  \sqrt{ \mathrm{e} }  {\mathbb E}_{|\theta} \left \{   \left ( 1 -  \mathrm{e}^{-t^{2} ( T_{2}^{\prime2} (\theta , n) -V_n^2 ( \theta )) /2 }\right )   I_{A_n \cap B_3(\theta)} \right \}  \, . \end{multline*}
 Again, using the previous arguments, we get that  
 \begin{equation}
 {\mathbb E}_{|\theta} ( Q^{\prime}_n(\theta)  I_{B_3(\theta)\cap A_n} )  \ll t^2 {\mathbb E}_{|\theta} \Big | \sum_{k=1}^n   \theta_k^2 \big ( {\mathbb E}_{k-1} (d_k^2) - {\mathbb E} (d_k^2)   \big ) \Big  |  \, . \label{estimate B3}%
\end{equation}
So, overall, taking into account \eqref{estimate B1}-\eqref{estimate B3}, we get that 
\begin{equation} \label{BoundIIn}
I\!\!I_{n}(\theta) \ll  t^2 {\mathbb E}_{|\theta} \Big | \sum_{k=1}^n   \theta_k^2 \big ( {\mathbb E}_{k-1} (d_k^2) - {\mathbb E} (d_k^2)   \big ) \Big  |  \, ,
\end{equation}
which combined with \eqref{BoundIn} entails that  
\begin{equation}
\big |{\mathbb E}_{|\theta} \big (  \Delta_{n}(\theta)I_{H_{n,K}^{\prime }(\theta) \cap A_n} \big ) \big |\ll |t|^{3}{\mathbb E}_{|\theta}(|T_{3}^{3}|)+t^{4}{\mathbb E}_{|\theta}(T_{4}^{4}) + t^2 {\mathbb E}_{|\theta} \Big | \sum_{k=1}^n   \theta_k^2 \big ( {\mathbb E}_{k-1} (d_k^2) - {\mathbb E} (d_k^2)   \big ) \Big  |  \, .\label{second part 0}%
\end{equation}
Starting from \eqref{first part 2} and taking into account \eqref{truc1} (combined with  \eqref{estimate Anc}) and \eqref{second part 0}, it follows that 
\begin{equation} \label{second part}
\big |{\mathbb E}_{|\theta} \big (  \Delta_{n}(\theta)I_{H_{n,K}^{\prime }(\theta) } \big ) \big |\ll |t|^{3}{\mathbb E}_{|\theta}(|T_{3}^{3}|)+t^{4}{\mathbb E}_{|\theta}(T_{4}^{4}) + t^2 {\mathbb E}_{|\theta}\Big | \sum_{k=1}^n   \theta_k^2 \big ( \E_{k-1} (d_k^2) - \E (d_k^2)   \big ) \Big  |  \, .
\end{equation}
Overall, starting from (\ref{decopose}) and combining \eqref{first part} and \eqref{second part}, the first part of the Lemma \ref{Lma1auxiliarymds} follows.

To prove the second part we use the following decomposition (see for instance (2.32) in \cite{MPU19}): 
\[
\left\vert {\mathbb E}_{|\theta} \left(  \mathrm{e}^{itS_{n}(\theta)}-\mathrm{e}^{-t^{2}V_n^2 ( \theta)/2}\right)
\right\vert \leq2{\mathbb P}_{|\theta}(A_{n}^{c})+\left\vert {\mathbb E}_{|\theta} \left(  \mathrm{e}^{itS_{n}^{\prime
}(\theta)}-\mathrm{e}^{-t^{2}V_n^2 ( \theta)/2}\right)  \right\vert \, .
\]
Then the second part of Lemma \ref{Lma1auxiliarymds} follows from the first part of Lemma \ref{Lma1auxiliarymds} and the upper bound \eqref{estimate Anc}. $\square$

%

\subsection{Proof of Theorems \ref{Th martingale} and \ref{Th martingale2}.}

We start with the following lemma concerning the vector
$\theta$:

\begin{lemma}
\label{lemma Moments teta} Assume that $\theta$ is uniformly distributed on
the sphere $S^{n-1}$ of ${\mathbb R}^n$ and $\left(  a_{j}\right)  _{1\leq j\leq n}$ is a vector in ${\mathbb R}^n$. Then, there exists $c>0$ such that for any positive integer $n$, 
\[
{\mathbb E} \left(\sum_{j=1}^{n}\theta_{j}^{4}|a_{j}| \right ) \leq c\frac{1}{n}\left(  \frac
{1}{n}\sum_{j=1}^{n}|a_{j}|\right)  ,
\]
and%
\[
{\mathbb E} \left\vert \sum_{j=1}^{n}\theta_{j}^{3}a_{j}\right\vert \leq
c\frac{1}{n}\left(  \frac{1}{n}\sum_{j=1}^{n}a_{j}^{2}\right)  ^{1/2}.
\]

\end{lemma}

\noindent{\bf Proof.} Let us introduce $Z=(Z_{1},\dots,Z_{n})$ i.i.d. standard normal variables
independent of $\theta.$ Let $\gamma$ be a chi-square random variable with $n$
degrees of freedom independent of $\theta$. Observe that the vector $Z$ has
the same distribution as $\sqrt{\gamma}\theta$. Since ${\mathbb E}(Z_{1})=0,$ by the
independence of the $Z_i$'s,
\[
{\mathbb E}\left\vert \sum_{j=1}^{n}Z_{j}^{3}a_{j}\right\vert ={\mathbb E} \left\vert \sum
_{j=1}^{n}\gamma^{3/2}\theta_{j}^{3}a_{j}\right\vert ={\mathbb E}(\gamma^{3/2}){\mathbb E} \left\vert \sum_{j=1}^{n}\theta_{j}^{3}a_{j}\right\vert .
\]
So%
\begin{align*}
{\mathbb E} \left\vert \sum_{j=1}^{n}\theta_{j}^{3}a_{j}\right\vert  &  \ll
\frac{1}{n^{3/2}}{\mathbb E}\left\vert \sum_{j=1}^{n}Z_{j}^{3}a_{j}\right\vert \leq
\frac{1}{n^{3/2}}\left \| \sum_{j=1}^{n}Z_{j}^{3}a_{j}\right \|_2
\\
&  \ll\frac{1}{n^{3/2}}\left(  \sum_{j=1}^{n}a_{j}^{2}\right)  ^{1/2} \, ,
\end{align*}
and the second inequality is proved.

Also, we can use  the distribution of $\theta_{j}$ which implies that there exists $C >0$ such that $\sup_{1 \leq i \leq n}  {\mathbb E} ( \theta_i^4) \leq C n^{-2}$. This implies directly the first part of the lemma. $  \square$

\bigskip

We combine now Lemma \ref{Lma1auxiliarymds} with Lemma
\ref{lemma Moments teta} and use the fact that ${\mathbb E} (\theta_k^2)= 1/n$ for any $1 \leq k \leq n$.  This gives Proposition \ref{prop dist Ch funct} below, that will be used to prove 
Theorem \ref{Th martingale}. 

\begin{proposition}
\label{prop dist Ch funct}  For any real $t$, we have
\begin{multline*}
{\mathbb E}\left\vert {\mathbb E}_{|\theta}\left(  \mathrm{e}^{itS_{n}(\theta)}-\mathrm{e}%
^{-t^{2}V^2_n(\theta)/2}\right)  \right\vert \ll\frac{1}{n}|t|^{3} {\mathbb E} \left(  \frac{1}{n}%
\sum_{j=1}^{n}\left ({\mathbb E}_{j-1} ( d_{j}^{3} ) \right)^2 \right)^{1/2}+t^{4}\frac{1}{n}\left(  \frac{1}%
{n}\sum_{j=1}^{n}{\mathbb E}(d_{j}^{4})\right)   \\ + t^2  \frac{1}{n} \sum_{k=1}^n \Vert {\mathbb E}_{k-1}(d_k^2) - {\mathbb E}(d_k^2) \Vert_1\,  .
\end{multline*}
\end{proposition}

\begin{remark} If ${\mathbb E}_{j-1} ( d_{j}^{2} ) = {\mathbb E}  ( d_{j}^{2} ) =1$, then 
\[
\left(  {\mathbb E}_{j-1}(d_{j}^{3})\right)  ^{2}\leq\left(  {\mathbb E}_{j-1}(d_{j}^{2})\right)
\left(  {\mathbb E}_{j-1}(d_{j}^{4})\right)  \leq {\mathbb E}_{j-1}(d_{j}^{4}) \,  ,
\]
implying that 
\[
{\mathbb E} \left(  \frac{1}{n}%
\sum_{j=1}^{n}\left ({\mathbb E}_{j-1} ( d_{j}^{3}  )  \right)^2 \right)  ^{1/2} \leq   \left(  \frac{1}{n}%
\sum_{j=1}^{n}{\mathbb E} ( d_{j}^{4} ) \right)  ^{1/2} \, .
\]
\end{remark}

\begin{remark} \label{conditioncondexpect}
Using the notation \eqref{Bobkov1}, 
it follows from Proposition \ref{prop dist Ch funct} that 
\[
{\mathbb E}\left\vert {\mathbb E}_{|\theta}\left(  \mathrm{e}^{itS_{n}(\theta)}-\mathrm{e}%
^{-t^{2}V^2_n(\theta)/2}\right)  \right\vert \ll\frac{1}{n}|t|^{3}  (  \alpha_n |t| + \beta_n)    + t^2  \frac{1}{n} \sum_{k=1}^n \Vert {\mathbb E}_{k-1}(d_k^2) - {\mathbb E}(d_k^2) \Vert_1\,  .
\]
\end{remark}

\bigskip 
For the proof of Theorem  \ref{Th martingale2}, we shall rather use Proposition \ref{propalter dist Ch funct} below. 

\begin{proposition}
\label{propalter dist Ch funct}  We have
\begin{multline*}
{\mathbb E}\left\vert {\mathbb E}_{|\theta} \left(  \mathrm{e}^{itS_{n}(\theta)}-\mathrm{e}%
^{-t^{2}V^2_n(\theta)/2}\right)  \right\vert \ll\frac{1}{n}|t|^{3} {\mathbb E} \left(  \frac{1}{n}%
\sum_{j=1}^{n}\left({\mathbb E}_{j-1} ( d_{j}^{3} ) \right )^2\right)^{1/2}+t^{4}\frac{1}{n}\left(  \frac{1}%
{n}\sum_{j=1}^{n}{\mathbb E}(d_{j}^{4})\right)   \\ +\frac{ t^2 }{\sqrt n}  \left(  \frac{1}%
{n}\sum_{j=1}^{n}{\mathbb E}(d_{j}^{4})\right)^{1/2} + \frac{ t^2 }{\sqrt n}  c_2(n) \,  .
\end{multline*}
where $c_2(n) $ is defined in \eqref{Bobkov3alter}
\end{proposition}
\noindent {\bf Proof.} We start again from Lemma \ref{Lma1auxiliarymds} combined with Lemma
\ref{lemma Moments teta}. Here,  compared to the statement of Proposition \ref{prop dist Ch funct}, we only have to show that 
\begin{equation}  \label{alterlma}
 \Big \Vert \sum_{k=1}^n   \theta_k^2 \big ( \E_{k-1} (d_k^2) - \E (d_k^2)   \big ) \Big  \Vert_1 \ll n^{-1/2}  \Big \{ \Big (  \frac{1}%
{n}\sum_{j=1}^{n}\E (d_{j}^{4})\Big )^{1/2} +  c_2(n) \Big \} \, .
\end{equation}
With this aim, we first write
\[
  \Big \Vert \sum_{k=1}^n   \theta_k^2 \big ( \E_{k-1} (d_k^2) - \E (d_k^2)   \big ) \Big  \Vert_1 \leq \Big \Vert \sum_{k=1}^n   \theta_k^2 \big ( \E_{k-1} (d_k^2) - \E (d_k^2)   \big ) \Big  \Vert_2  \, .
\]
Next,  
\begin{multline*}
{\mathbb E}_{|\theta} \Big ( \sum_{k=1}^n   \theta_k^2 \big ( \E_{k-1} (d_k^2) - \E (d_k^2)   \big ) \Big  )^2 = {\rm Var }_{|\theta}  \Big ( \sum_{k=1}^n   \theta_k^2 \E_{k-1} (d_k^2)  \Big )  \\
 \leq 2 {\rm Var }_{|\theta}  \Big ( \sum_{k=1}^n   \theta_k^2  \big ( \E_{k-1} (d_k^2)  -d_k^2  \big )   \Big ) + 2  {\rm Var }_{|\theta}  \Big ( \sum_{k=1}^n   \theta_k^2 d_k^2  \Big ) \, .
\end{multline*}
Now
\[
 {\rm Var }_{|\theta}  \Big ( \sum_{k=1}^n   \theta_k^2  \big (  d_k^2 - \E_{k-1} (d_k^2)  \big )   \Big ) =  \sum_{k=1}^n    \theta_k^4   {\rm Var }_{|\theta}  (    d_k^2 - \E_{k-1} (d_k^2))  
 \leq 4   \sum_{k=1}^n    \theta_k^4 \E (d_k^4) \, . 
\]
So, overall, 
\[
 {\mathbb E}_{|\theta} \Big (\sum_{k=1}^n   \theta_k^2 \big ( \E_{k-1} (d_k^2) - \E (d_k^2)   \big ) \Big  )^2 \leq   4   \sum_{k=1}^n    \theta_k^4 \E (d_k^4) + 
 2   \sum_{k, \ell=1}^n   \theta_k^4   | {\rm cov} ( d_k^2, d_\ell^2) |  \, .
\]
Next we integrate with respect to $\theta$. Using the fact that, for any $1 \leq k \leq n$,  $\E ( \theta_k^4) \ll n^{-2}$, inequality 
\eqref{alterlma}  follows. $\square$

\bigskip

Now we have to relate Propositions \ref{prop dist Ch funct} and \ref{propalter dist Ch funct} to the Kolmogorov distance. With this aim, we use the approach in Bobkov-Chistyakov-G\"otze  \cite{BCG20b} (see their Lemma 6.1).  

We have the following proposition:

\begin{proposition}
\label{smoothingLemma} Let $(d_{k})_{1\leq k\leq n}$ be a vector of martingale
differences in ${\mathbb L}^4$ and $N_\theta$ a centered  normal variable with variance $V_n^2(\theta) = \sum_{i=1}^n \theta_i^2 \E (d_i^2)$. If  
(\ref{Bobkov2}) is satisfied then, for all $T\geq
T_{0}>1$,
\begin{multline}
{\E}\left(  \kappa_{\theta}(P_{S_{n}(\theta)}, P_{N_\theta})\right)  \ll\int
_{0}^{T_{0}}{\E}|\E_{|\theta}(\mathrm{e}^{itS_{n}(\theta)})-\mathrm{e}^{-t^{2}V_n^2(\theta)  
/2}|\frac{dt}{t} \\
+ ( \alpha_n + \sigma_{4}^{2} (n))  \frac{ (1+\log(T/T_{0}))}{n}+\mathrm{e}^{-T_{0}^{2}/16}+\frac
{1}{T}\, . \label{B2}%
\end{multline}
\end{proposition}

\noindent{\bf Proof.} By the usual Berry-Esseen smoothing inequality (see for instance inequality (6.2) in \cite{BCG18}), there exists  $c>0$ such that, for any $T>0$, 
\begin{equation} \label{BES}
c \E \Big (  \kappa_{\theta}(P_{S_{n}(\theta)},P_{N_\theta} )  \Big ) \leq \int_0^T \E \big |  \E_{| \theta}(\mathrm{e}^{itS_{n}(\theta)}) - n_{\theta} (t) \big |  \frac{dt }{t } + \frac{1}{T} \int_0^T \E  |n_{\theta} (t) | dt, 
\end{equation}
where $n_{\theta} (t) = {\mathbb E}_{|\theta} ( \mathrm{e}^{it N_\theta} ) =  \mathrm{e}^{-t^2 \sum_{i=1}^n \theta_i^2 \E (d_i^2)/2}=  \mathrm{e}^{-t^2 V_n^2(\theta) /2} $.  Since  (\ref{Bobkov2}) is assumed,  Lemma 5.2 together with Corollary 2.3 in \cite{BCG18} imply that for
all $T\geq T_{0}>1$,
\[
\int_{T_{0}}^{T}t^{-1}{\E}|\E_{|\theta}(\mathrm{e}^{itS_{n}(\theta)}%
)|dt\ll ( M^4_{4}(n)  + \sigma_{4}^{2} (n))   \frac{ (1+ \log(T/T_{0}))}{n}+\mathrm{e}^{-T_{0}^{2}/16}\,, \label{B1}%
\]
where 
\[
M^4_{4}(n) =\left \|{\E_{|\theta}}\left ( \Big (\sum_{k=1}^{n}\theta_{k}d_{k} \Big )^{4} \right ) \right \|_\infty \, .
\]
For martingales, by the Burkholder's inequality as given by Rio \cite{Rio09} (see also \cite{MPU19}, Theorem 2.13),  we have
\[
\E_{|\theta}\left (\Big (  \sum_{k=1}^{n}\theta_{k}d_{k} \Big )^{4}\right )\leq 3^2 \left(  \sum_{k=1}%
^{n}\left (\E_{|\theta}  (\theta_{k}^4d_{k}^{4} ) \right )^{1/2}\right)  ^{2}\leq 9 \sup_{1 \leq k \leq n}\E(d_{k}^{4}%
) = 9 \alpha_n  \, .
\]
Therefore, if    (\ref{Bobkov2}) holds,  we have, for
all $T\geq T_{0}>1$,
\begin{equation}
\int_{T_{0}}^{T}t^{-1}{\E}|\E_{|\theta}(\mathrm{e}^{itS_{n}(\theta)}%
)|dt\ll ( \alpha_n + \sigma_{4}^{2} (n))   \frac{ (1+ \log(T/T_{0}))}{n}+\mathrm{e}^{-T_{0}^{2}/16}\, . \label{B1}%
\end{equation}
Next, note that we can also write $n_{\theta} (t) = \E_{|\theta}(\mathrm{e}^{itU_{n}(\theta)}%
)$ where $U_{n}(\theta)= \sum_{i=1}^n \theta_i Y_i$ and the $Y_i$'s are independent and such that $Y_i \sim {\mathcal N}(0, \E (d_i^2) )$.  Now
\[
\sup_{\theta} \E_{|\theta}  (U_{n}(\theta)^4) = 3  \sup_{\theta}  \left ( \E_{|\theta}  (U_{n}(\theta)^2)\right )^2  \leq 3   \left (  \sup_{1 \leq k \leq n}\E(d_{k}^{2}) \right ) ^2 \leq  3 \alpha_n \, , 
\]
and 
\[
n^{-1} \mathrm{Var}\left (\sum_{k=1}^{n}Y_{k}^{2} \right )  = n^{-1} \sum_{k=1}^{n} \mathrm{Var} ( Y_k^2) \leq 3  n^{-1} \sum_{k=1}^{n} (  \E (Y_k^2) )^2  \leq 3 \alpha_n \, .
\]
Therefore, taking into account (\ref{Bobkov2}),   Lemma 5.2 together with Corollary 2.3 in \cite{BCG18} imply that
\begin{equation*} \label{Boundntheta}
 |n_{\theta} (t) |  \ll  n^{-1 } \alpha_n + e^{-t^2/16} \, ,
\end{equation*}
and consequently,  for
all $T\geq T_{0}>1$, 
\begin{equation}
\int_{T_{0}}^{T}t^{-1}{\E}|n_{\theta} (t)|dt + T^{-1} \int_0^T {\E}|n_{\theta} (t)|dt   \ll  \alpha_n   \frac{ (1+ \log(T/T_{0}))}{n}+\mathrm{e}^{-T_{0}^{2}/16}
+ \frac{1}{T}\, . \label{B1G}%
\end{equation}
Starting from \eqref{BES} and taking into account \eqref{B1} and \eqref{B1G}, the proposition follows.  $\square$

\subsection{End of the proof of Theorems \ref{Th martingale} and \ref{Th martingale2}.} We shall apply Proposition \ref{smoothingLemma} with $T=n$ and $T_{0}%
= 4 \sqrt{\log n}$. Therefore, to derive an upper bound for ${\E}\left (
\kappa_{\theta}(P_{S_{n}(\theta)},P_{N_\theta} )  \right )  $ of order $1/n$ modulo an
extra term $(\log n)^{2}$, one needs to prove that
\[
\int_{0}^{T_{0}}{\E}|\E_{|\theta}(\mathrm{e}^{itS_{n}(\theta)})-\mathrm{e}%
^{-t^{2}V_n^2 (\theta) /2}|\frac{dt}{t}\ll \frac{1}{n}\left( \gamma_n (\log n )+  \beta_n (\log n)^{3/2}+ \alpha_n(\log n)^{2}\right)   \,.
\]
By  Proposition \ref{prop dist Ch funct},
\begin{align*}
\int_{0}^{T_{0}}{\E}|\E_{|\theta}(\mathrm{e}^{itS_{n}(\theta)})-\mathrm{e}%
^{-t^{2}V_n^2 (\theta) / 2}|\frac{dt}{t}  &  \ll\frac{1}{n}\int_{0}^{T_{0}}(t^{2} \beta_n+t^{3} \alpha_n
 + t \gamma_n )dt\ll\frac{1}{n}\left(  |T_{0}|^{3}+|T_{0}|^{4}\right) \\
&  \ll\frac{1}{n}\left( \gamma_n (\log n )+  \beta_n (\log n)^{3/2}+ \alpha_n(\log n)^{2}\right) \, ,
\end{align*}
and  Theorem \ref{Th martingale} follows. $\square$

The proof of Theorem  \ref{Th martingale2} is similar to the one of Theorem \ref{Th martingale} above. The only difference is that we use Proposition \ref{propalter dist Ch funct} instead of Proposition \ref{prop dist Ch funct} and we take into account Remark \ref{remarkonGaussian}.  $\square$

\subsection{Proof of Remark \ref{remarkonGaussian}} \label{Sectionremark2} We start by proving \eqref{boundGaussian}.  With this aim, we  recall Proposition 2.6 in Johnston-Prochno \cite{JP22}. Let $G_{\sigma}$ be a centered Gaussian r.v. with variance $\sigma^2$. Let $\alpha, \beta $ in $(0, \infty)$. If $\beta / \alpha > 1/2$, then 
\begin{equation} \label{ineqJP0}
\sup_{t \in {\mathbb R}} \big | {\mathbb P} ( G_{\alpha } \leq t)- {\mathbb P} ( G_{\beta } \leq t) \big | \leq \frac{3}{8} \frac{| \alpha^2 - \beta^2|}{\alpha^2} \, .
\end{equation}
Denote $V_{n}^{2}(\theta) =\theta_{1}^{2}b_{1}^{2}+\cdots+\theta_{n}^{2}b_{n}^{2}$
with $b_{j}^{2}=\E(d_{j}^{2})$ for all $j$.   Using inequality \eqref{ineqJP0} with $\alpha =1$ and $\beta^2= V_{n}^{2}(\theta)$, we get 
\begin{multline*}
{\E} \big (
\kappa_{\theta}(P_{N_\theta},P_{N} ) \big ) = {\E}  \big (
\kappa_{\theta}(P_{N_\theta},P_{N} )  {\bf 1}_{\{ V_{n}(\theta) > 1/2  \} }\big ) + {\E}  \big (
\kappa_{\theta}(P_{N_\theta},P_{N} ) {\bf 1}_{\{ V_{n}(\theta)  \leq 1/2  \} }\big )   \\
\leq  \frac{3}{8} \E \big |  V_{n}^{2}(\theta)  - 1 \big |  +  {\mathbb P}( V_{n}(\theta)  \leq 1/2) \, .
\end{multline*}
Now
\[
{\mathbb P} ( V_{n}(\theta)  \leq 1/2) = {\mathbb P} (  1 -  V^2_{n}(\theta)  \geq 1/4) \leq  {\mathbb P} ( | 1 -  V^2_{n}(\theta)  | \geq 1/4) \leq 4 \E \big |  V_{n}^{2}(\theta)  - 1\big | \, .
\]
Hence
\begin{equation} \label{ineqJP}
{\E}\big (
\kappa_{\theta}(P_{N_\theta},P_{N} ) \big )  \leq \frac{35}{8} \E \big |  V_{n}^{2}(\theta)  - 1 \big | \, .
\end{equation}
%
%
To give an upper bound of the term in the right-hand side, we are going to use a second order concentration on the sphere developed in
\cite{BCG17} and also given in Proposition 4.1 in \cite{BCG20a}. Let 
\[
v(\theta)=V_{n}^{2}(\theta)-1=\sum_{k=1}^{n}\theta_{k}^{2}b_{k}^{2}-1.
\]
Notice that  $v(\theta)$ has mean $0$ because $\sum_{k=1}^{n}b_{k}^{2}=n$ and
$\E (\theta_{k}^{2})=1/n$ for all $k$. Moreover it is an even function 
($v(\theta)=v(-\theta)$), so it is orthogonal  in the Hilbert space ${\mathbb L}^2 ( \mu_{n-1} )$ to all linear functions on the sphere. In addition, it is continuous
and infinitely  differentiable for all $\theta$ and the derivatives are
continuous on an open interval containing $\theta.$ Proposition 4.1 in \cite{BCG20a} then asserts that, for any
$a$ in ${\mathbb R}$,  
\begin{equation}
\E \left (|v(\theta)|^{2} \right) \leq\frac{5}{\left(
n-1\right)  ^{2}}\E \left (\|\nabla^{2}v(\theta)-aI_{n}\|_{HS}%
^{2}\right ),\label{PI}%
\end{equation}
where $HS$ stands for the Hilbert-Schmidt norm and $I_n$ denotes the identity matrix of size $n$. We have 
\[
\frac{\partial^{2}v(\theta)}{\partial\theta_{k}\partial\theta_{j}}%
=2b_{k}^{2}\delta_{j,k},
\]
with $\delta_{j,k}=1,$ for $j=k$ and $0$ in rest. So, taking $a=2$, 
\[
\E \left (\|\nabla^{2}v (\theta)-aI_{n}\|_{HS}%
^{2}\right )= 4 \sum_{k=1}^{n} (  b_{k}^{2}- 1 )^{2} \, .
\]
This implies that, for any $n \geq 2$, 
\begin{equation} \label{R0}
 \left \|\sum_{k=1}^{n}\theta_{k}^{2}b_{k}^{2}-1 \right \|_1 \leq  \left \|\sum_{k=1}^{n}\theta_{k}^{2}b_{k}^{2}-1\right \|_2\leq 2 \sqrt{5} \frac{ \Big (  \sum_{k=1}^{n} (  b_{k}^{2}- 1  )^2 \Big )^{1/2}   }{n-1} \, .
\end{equation}
Starting from \eqref{ineqJP} and considering the upper bound above, inequality \eqref{boundGaussian} follows. 


To prove \eqref{boundGaussian2}, we start with the usual Berry-Esseen smoothing inequality which implies that there exists a $c>0$ such that, for any $T>0$, 
\begin{equation} \label{BESalter}
c  \sup_{t\in{\mathbb{R}}}| {\mathbb P}(N_\theta\leq t)-{\mathbb P}(N \leq t)|  \leq \int_0^T  \big |  \E(\mathrm{e}^{itU_{n}(\theta)}) - \mathrm{e}%
^{-t^{2}/2}  \big |  \frac{dt }{t } + \frac{1}{T} \int_0^T  \mathrm{e}%
^{-t^{2}/2}  dt  \, , 
\end{equation}
where $U_{n}(\theta)= \sum_{i=1}^n \theta_i Y_i$ and the $Y_i$'s are independent, independent of $\theta$ and such that $Y_i \sim {\mathcal N}(0, \E (d_i^2) )$.  

Taking into account \eqref{B1} with $S_{n}(\theta)$ replaced by $U_{n}(\theta)$ and setting $T_0 = 4 \sqrt{\log n}$, we infer that the upper bound  \eqref{boundGaussian2} 
follows from \eqref{BESalter} if one can prove  that 
\begin{equation} \label{R2}
\int_{0}^{T_{0}} \left | {\E} \left (\mathrm{e}%
^{-t^{2}V_n^2 (\theta) /2}\right ) - \mathrm{e}%
^{-t^{2}/2} \right |\frac{dt}{t}  \ll  \frac{ \alpha_n ( \log n)^2   }{n}\, .
\end{equation}
With this aim note that 
\begin{multline*}
{\E}  \left (\mathrm{e}%
^{-t^{2}V_n^2 (\theta) /2} \right )
 -  {\rm e}^{-t^2    /2}    
 = - \frac{t^2}{2}  {\E}  \big (  V_n^2(\theta)-  1 \big )     {\rm e}^{-t^2      /2}  \\+  \frac{t^4}{8 }   {\E} \Big [     \big (  V_n^2(\theta)-  1 \big )^2  \int_0^1 (1-s)  
   {\rm e}^{-t^2   [ 1  +s (  V_n^2(\theta) - 1)  ]   /2}  ds  \Big ] \, .
\end{multline*}
Now ${\E}  \left (  V_n^2(\theta) \right ) =  \sum_{i=1}^n \E ( \theta_i^2) \E (d_i^2) = n^{-1}   \sum_{i=1}^n  \E (d_i^2) =1   $. Hence
\[
\Big | 
 {\E}\left (   {\rm e}^{-t^2    V_n^2 ( \theta)   /2}     \right )  
 -  {\rm e}^{-t^2    /2}     \Big |  \leq   \frac{t^4}{8 }   {\E} \left [     \big (  V_n^2(\theta)- 1 \big )^2   \right ]   \, ,
\]
which combined with \eqref{R0} implies that, for any $n \geq 2$, 
\[
\Big | 
 {\E} \left (   {\rm e}^{-t^2    V_n^2 ( \theta)   /2}     \right )  
 -  {\rm e}^{-t^2    /2}     \Big |  \leq   \frac{5 t^4}{2 }  \frac{ \sum_{k=1}^{n} (  b_{k}^{2}- 1  )^2   }{(n-1)^2}  \ll   t^4  \frac{ \alpha_n   }{n}   \, . 
\]
This ends the proof of \eqref{R2}.  $\square$

\end{document}